\documentclass[12pt]{article}
\usepackage{amsmath,amssymb,amsthm,latexsym,amscd}

\title{Ranks for families of theories and their spectra\footnote{
\newline\indent \ \ \ This research was partially
supported by Russian Foundation for Basic Researches (Project No.
17-01-00531-a) and Committee of Science in Education and Science
Ministry of the Republic of Kazakhstan (Grant No. AP05132546).} }
\author{Sergey V.
Sudoplatov\footnote{sudoplat@math.nsc.ru}}
\date{}
\begin{document}
\maketitle

\begin{abstract}
We define ranks and degrees for families of theories, similar to
Morley rank and degree, as well as Cantor-Bendixson rank and
degree, and the notion of totally transcendental family of
theories. Bounds for $e$-spectra with respect to ranks and degrees
are found. It is shown that the ranks and the degrees are
preserved under $E$-closures and values for the ranks and the
degrees are characterized. Criteria for totally transcendental
families in terms of cardinality of $E$-closure and of the
$e$-spectrum value, for a countable language, are proved.

\end{abstract}

We continue to study families of theories \cite{cs, cl, lut, ccct,
rest, lft} and their approximations \cite{at} introducing ranks
and degrees for families of theories, similar to Morley rank and
degree \cite{Morley}, as well as Cantor-Bendixson rank and degree,
and the notion of totally transcendental family of theories. These
ranks and degree plays a similar role for families of theories,
with hierarchies for definable sets of theories, as Morley ones
for a fixed theory although they have own specificities.

Bounds for $e$-spectra with respect to ranks and degrees are
found. It is shown that the ranks and the degrees are preserved
under $E$-closures and values for the ranks and the degrees are
characterized. Criteria for totally transcendental families in
terms of cardinality of $E$-closure and of the $e$-spectrum value,
for a countable language, are proved.

\section{Preliminaries}

Throughout the paper we consider complete first-order theories $T$
in predicate languages $\Sigma(T)$ and use the following
terminology in \cite{cs, cl, lut, ccct, rest, lft}.

Let $P=(P_i)_{i\in I}$, be a family of nonempty unary predicates,
$(\mathcal{A}_i)_{i\in I}$ be a family of structures such that
$P_i$ is the universe of $\mathcal{A}_i$, $i\in I$, and the
symbols $P_i$ are disjoint with languages for the structures
$\mathcal{A}_j$, $j\in I$. The structure
$\mathcal{A}_P\rightleftharpoons\bigcup\limits_{i\in
I}\mathcal{A}_i$\index{$\mathcal{A}_P$} expanded by the predicates
$P_i$ is the {\em $P$-union}\index{$P$-union} of the structures
$\mathcal{A}_i$, and the operator mapping $(\mathcal{A}_i)_{i\in
I}$ to $\mathcal{A}_P$ is the {\em
$P$-operator}\index{$P$-operator}. The structure $\mathcal{A}_P$
is called the {\em $P$-combination}\index{$P$-combination} of the
structures $\mathcal{A}_i$ and denoted by ${\rm
Comb}_P(\mathcal{A}_i)_{i\in I}$\index{${\rm
Comb}_P(\mathcal{A}_i)_{i\in I}$} if
$\mathcal{A}_i=(\mathcal{A}_P\upharpoonright
A_i)\upharpoonright\Sigma(\mathcal{A}_i)$, $i\in I$. Structures
$\mathcal{A}'$, which are elementary equivalent to ${\rm
Comb}_P(\mathcal{A}_i)_{i\in I}$, will be also considered as
$P$-combinations.

Clearly, all structures $\mathcal{A}'\equiv {\rm
Comb}_P(\mathcal{A}_i)_{i\in I}$ are represented as unions of
their restrictions $\mathcal{A}'_i=(\mathcal{A}'\upharpoonright
P_i)\upharpoonright\Sigma(\mathcal{A}_i)$ if and only if the set
$p_\infty(x)=\{\neg P_i(x)\mid i\in I\}$ is inconsistent. If
$\mathcal{A}'\ne{\rm Comb}_P(\mathcal{A}'_i)_{i\in I}$, we write
$\mathcal{A}'={\rm Comb}_P(\mathcal{A}'_i)_{i\in
I\cup\{\infty\}}$, where
$\mathcal{A}'_\infty=\mathcal{A}'\upharpoonright
\bigcap\limits_{i\in I}\overline{P_i}$, maybe applying
Morleyzation. Moreover, we write ${\rm
Comb}_P(\mathcal{A}_i)_{i\in I\cup\{\infty\}}$\index{${\rm
Comb}_P(\mathcal{A}_i)_{i\in I\cup\{\infty\}}$} for ${\rm
Comb}_P(\mathcal{A}_i)_{i\in I}$ with the empty structure
$\mathcal{A}_\infty$.

Note that if all predicates $P_i$ are disjoint, a structure
$\mathcal{A}_P$ is a $P$-combination and a disjoint union of
structures $\mathcal{A}_i$. In this case the $P$-combination
$\mathcal{A}_P$ is called {\em
disjoint}.\index{$P$-combination!disjoint} Clearly, for any
disjoint $P$-combination $\mathcal{A}_P$, ${\rm
Th}(\mathcal{A}_P)={\rm Th}(\mathcal{A}'_P)$, where
$\mathcal{A}'_P$ is obtained from $\mathcal{A}_P$ replacing
$\mathcal{A}_i$ by pairwise disjoint
$\mathcal{A}'_i\equiv\mathcal{A}_i$, $i\in I$. Thus, in this case,
similar to structures the $P$-operator works for the theories
$T_i={\rm Th}(\mathcal{A}_i)$ producing the theory $T_P={\rm
Th}(\mathcal{A}_P)$\index{$T_P$}, being {\em
$P$-combination}\index{$P$-combination} of $T_i$, which is denoted
by ${\rm Comb}_P(T_i)_{i\in I}$.\index{${\rm Comb}_P(T_i)_{i\in
I}$}

Notice that $P$-combinations are represented by generalized
products of structures \cite{FV}.

For an equivalence relation $E$ replacing disjoint predicates
$P_i$ by $E$-classes we get the structure
$\mathcal{A}_E$\index{$\mathcal{A}_E$} being the {\em
$E$-union}\index{$E$-union} of the structures $\mathcal{A}_i$. In
this case the operator mapping $(\mathcal{A}_i)_{i\in I}$ to
$\mathcal{A}_E$ is the {\em $E$-operator}\index{$E$-operator}. The
structure $\mathcal{A}_E$ is also called the {\em
$E$-combination}\index{$E$-combination} of the structures
$\mathcal{A}_i$ and denoted by ${\rm Comb}_E(\mathcal{A}_i)_{i\in
I}$\index{${\rm Comb}_E(\mathcal{A}_i)_{i\in I}$}; here
$\mathcal{A}_i=(\mathcal{A}_E\upharpoonright
A_i)\upharpoonright\Sigma(\mathcal{A}_i)$, $i\in I$. Similar
above, structures $\mathcal{A}'$, which are elementary equivalent
to $\mathcal{A}_E$, are denoted by ${\rm
Comb}_E(\mathcal{A}'_j)_{j\in J}$, where $\mathcal{A}'_j$ are
restrictions of $\mathcal{A}'$ to its $E$-classes. The
$E$-operator works for the theories $T_i={\rm Th}(\mathcal{A}_i)$
producing the theory $T_E={\rm Th}(\mathcal{A}_E)$\index{$T_E$},
being {\em $E$-combination}\index{$E$-combination} of $T_i$, which
is denoted by ${\rm Comb}_E(T_i)_{i\in I}$\index{${\rm
Comb}_E(T_i)_{i\in I}$} or by ${\rm
Comb}_E(\mathcal{T})$\index{${\rm Comb}_E(\mathcal{T})$}, where
$\mathcal{T}=\{T_i\mid i\in I\}$.

Clearly, $\mathcal{A}'\equiv\mathcal{A}_P$ realizing $p_\infty(x)$
is not elementary embeddable into $\mathcal{A}_P$ and can not be
represented as a disjoint $P$-combination of
$\mathcal{A}'_i\equiv\mathcal{A}_i$, $i\in I$. At the same time,
there are $E$-combinations such that all
$\mathcal{A}'\equiv\mathcal{A}_E$ can be represented as
$E$-combinations of some $\mathcal{A}'_j\equiv\mathcal{A}_i$. We
call this representability of $\mathcal{A}'$ to be the {\em
$E$-representability}.

If there is $\mathcal{A}'\equiv\mathcal{A}_E$ which is not
$E$-representable, we have the $E'$-representability replacing $E$
by $E'$ such that $E'$ is obtained from $E$ adding equivalence
classes with models for all theories $T$, where $T$ is a theory of
a restriction $\mathcal{B}$ of a structure
$\mathcal{A}'\equiv\mathcal{A}_E$ to some $E$-class and
$\mathcal{B}$ is not elementary equivalent to the structures
$\mathcal{A}_i$. The resulting structure $\mathcal{A}_{E'}$ (with
the $E'$-representability) is a {\em
$e$-completion}\index{$e$-completion}, or a {\em
$e$-saturation}\index{$e$-saturation}, of $\mathcal{A}_{E}$. The
structure $\mathcal{A}_{E'}$ itself is called {\em
$e$-complete}\index{Structure!$e$-complete}, or {\em
$e$-saturated}\index{Structure!$e$-saturated}, or {\em
$e$-universal}\index{Structure!$e$-universal}, or {\em
$e$-largest}\index{Structure!$e$-largest}.

For a structure $\mathcal{A}_E$ the number of {\em
new}\index{Structure!new} structures with respect to the
structures $\mathcal{A}_i$, i.~e., of the structures $\mathcal{B}$
which are pairwise elementary non-equivalent and elementary
non-equivalent to the structures $\mathcal{A}_i$, is called the
{\em $e$-spectrum}\index{$e$-spectrum} of $\mathcal{A}_E$ and
denoted by $e$-${\rm Sp}(\mathcal{A}_E)$.\index{$e$-${\rm
Sp}(\mathcal{A}_E)$} The value ${\rm sup}\{e$-${\rm
Sp}(\mathcal{A}'))\mid\mathcal{A}'\equiv\mathcal{A}_E\}$ is called
the {\em $e$-spectrum}\index{$e$-spectrum} of the theory ${\rm
Th}(\mathcal{A}_E)$ and denoted by $e$-${\rm Sp}({\rm
Th}(\mathcal{A}_E))$.\index{$e$-${\rm Sp}({\rm
Th}(\mathcal{A}_E))$} If structures $\mathcal{A}_i$ represent
theories $T_i$ of a family $\mathcal{T}$, consisting of $T_i$,
$i\in I$, then the $e$-spectrum $e$-${\rm Sp}(\mathcal{A}_E)$ is
denoted by $e$-${\rm Sp}(\mathcal{T})$.

If $\mathcal{A}_E$ does not have $E$-classes $\mathcal{A}_i$,
which can be removed, with all $E$-classes
$\mathcal{A}_j\equiv\mathcal{A}_i$, preserving the theory ${\rm
Th}(\mathcal{A}_E)$, then $\mathcal{A}_E$ is called {\em
$e$-prime}\index{Structure!$e$-prime}, or {\em
$e$-minimal}\index{Structure!$e$-minimal}.

For a structure $\mathcal{A}'\equiv\mathcal{A}_E$ we denote by
${\rm TH}(\mathcal{A}')$ the set of all theories ${\rm
Th}(\mathcal{A}_i)$\index{${\rm Th}(\mathcal{A}_i)$} of
$E$-classes $\mathcal{A}_i$ in $\mathcal{A}'$.

By the definition, an $e$-minimal structure $\mathcal{A}'$
consists of $E$-classes with a minimal set ${\rm
TH}(\mathcal{A}')$. If ${\rm TH}(\mathcal{A}')$ is the least for
models of ${\rm Th}(\mathcal{A}')$ then $\mathcal{A}'$ is called
{\em $e$-least}.\index{Structure!$e$-least}

\medskip
{\bf Definition} \cite{cl}. Let $\overline{\mathcal{T}}_\Sigma$ be
the set of all complete elementary theories of a relational
language $\Sigma$. For a set
$\mathcal{T}\subset\overline{\mathcal{T}}_\Sigma$ we denote by
${\rm Cl}_E(\mathcal{T})$ the set of all theories ${\rm
Th}(\mathcal{A})$, where $\mathcal{A}$ is a structure of some
$E$-class in $\mathcal{A}'\equiv\mathcal{A}_E$,
$\mathcal{A}_E={\rm Comb}_E(\mathcal{A}_i)_{i\in I}$, ${\rm
Th}(\mathcal{A}_i)\in\mathcal{T}$. As usual, if $\mathcal{T}={\rm
Cl}_E(\mathcal{T})$ then $\mathcal{T}$ is said to be {\em
$E$-closed}.\index{Set!$E$-closed}

The operator ${\rm Cl}_E$ of $E$-closure can be naturally extended
to the classes $\mathcal{T}\subset\overline{\mathcal{T}}$, where
$\overline{\mathcal{T}}$ is the union of all
$\overline{\mathcal{T}}_\Sigma$ as follows: ${\rm
Cl}_E(\mathcal{T})$ is the union of all ${\rm
Cl}_E(\mathcal{T}_0)$ for subsets
$\mathcal{T}_0\subseteq\mathcal{T}$, where new language symbols
with respect to the theories in $\mathcal{T}_0$ are empty.

For a set $\mathcal{T}\subset\overline{\mathcal{T}}$ of theories
in a language $\Sigma$ and for a sentence $\varphi$ with
$\Sigma(\varphi)\subseteq\Sigma$ we denote by
$\mathcal{T}_\varphi$\index{$\mathcal{T}_\varphi$} the set
$\{T\in\mathcal{T}\mid\varphi\in T\}$. Any set
$\mathcal{T}_\varphi$ is called the {\em $\varphi$-neighbourhood},
or simply a {\em neighbourhood}, for $\mathcal{T}$, or the
($\varphi$-){\em definable} subset of $\mathcal{T}$.

\medskip
{\bf Proposition 1.1} \cite{cl}. {\em If
$\mathcal{T}\subset\overline{\mathcal{T}}$ is an infinite set and
$T\in\overline{\mathcal{T}}\setminus\mathcal{T}$ then $T\in{\rm
Cl}_E(\mathcal{T})$ {\rm (}i.e., $T$ is an {\sl accumulation
point} for $\mathcal{T}$ with respect to $E$-closure ${\rm
Cl}_E${\rm )} if and only if for any formula $\varphi\in T$ the
set $\mathcal{T}_\varphi$ is infinite.}

\medskip
If $T$ is an accumulation point for $\mathcal{T}$ then we also say
that $T$ is an {\em accumulation point} for ${\rm
Cl}_E(\mathcal{T})$.

\medskip
{\bf Theorem 1.2} \cite{cl}. {\em For any sets
$\mathcal{T}_0,\mathcal{T}_1\subset\overline{\mathcal{T}}$, ${\rm
Cl}_E(\mathcal{T}_0\cup\mathcal{T}_1)={\rm
Cl}_E(\mathcal{T}_0)\cup{\rm Cl}_E(\mathcal{T}_1)$.}

\medskip
{\bf Definition} \cite{cl}. Let $\mathcal{T}_0$ be a closed set in
a topological space $(\mathcal{T},\mathcal{O}_E(\mathcal{T}))$,
where $\mathcal{O}_E(\mathcal{T})=\{\mathcal{T}\setminus{\rm
Cl}_E(\mathcal{T}')\mid\mathcal{T}'\subseteq\mathcal{T}\}$. A
subset $\mathcal{T}'_0\subseteq\mathcal{T}_0$ is said to be {\em
generating}\index{Set!generating} if $\mathcal{T}_0={\rm
Cl}_E(\mathcal{T}'_0)$. The generating set $\mathcal{T}'_0$ (for
$\mathcal{T}_0$) is {\em minimal}\index{Set!generating!minimal} if
$\mathcal{T}'_0$ does not contain proper generating subsets. A
minimal generating set $\mathcal{T}'_0$ is {\em
least}\index{Set!generating!least} if $\mathcal{T}'_0$ is
contained in each generating set for $\mathcal{T}_0$.

\medskip
{\bf Theorem 1.3} \cite{cl}. {\em If $\mathcal{T}'_0$ is a
generating set for a $E$-closed set $\mathcal{T}_0$ then the
following conditions are equivalent:

$(1)$ $\mathcal{T}'_0$ is the least generating set for
$\mathcal{T}_0$;

$(2)$ $\mathcal{T}'_0$ is a minimal generating set for
$\mathcal{T}_0$;

$(3)$ any theory in $\mathcal{T}'_0$ is isolated by some set
$(\mathcal{T}'_0)_\varphi$, i.e., for any $T\in\mathcal{T}'_0$
there is $\varphi\in T$ such that
$(\mathcal{T}'_0)_\varphi=\{T\}$;

$(4)$ any theory in $\mathcal{T}'_0$ is isolated by some set
$(\mathcal{T}_0)_\varphi$, i.e., for any $T\in\mathcal{T}'_0$
there is $\varphi\in T$ such that
$(\mathcal{T}_0)_\varphi=\{T\}$.}

\medskip
Notice that having the least generating set $\mathcal{T}'_0$ for a
$E$-closed set $\mathcal{T}_0$, $$e\mbox{-}{\rm
Sp}(\mathcal{T}_0)=e\mbox{-}{\rm
Sp}(\mathcal{T}'_0)=|\mathcal{T}_0\setminus\mathcal{T}'_0|.$$

\medskip
{\bf Definition} \cite{at}. Let $\mathcal{T}$ be a class of
theories and $T$ be a theory, $T\notin\mathcal{T}$. The theory $T$
is called {\em $\mathcal{T}$-approximated}, or {\em approximated
by} $\mathcal{T}$, or {\em $\mathcal{T}$-approximable}, or a {\em
pseudo-$\mathcal{T}$-theory}, if for any formula $\varphi\in T$
there is $T'\in\mathcal{T}$ such that $\varphi\in T'$.

If $T$ is $\mathcal{T}$-approximated then $\mathcal{T}$ is called
an {\em approximating family} for $T$, theories $T'\in\mathcal{T}$
are {\em approximations} for $T$, and $T$ is an {\em accumulation
point} for $\mathcal{T}$.

An approximating family $\mathcal{T}$ is called {\em
single-valued}, or {\em $e$-categorical}, if $e$-${\rm
Sp}(\mathcal{T})=1$.

An approximating family $\mathcal{T}$ is called {\em $e$-minimal}
if for any sentence $\varphi\in\Sigma(T)$, $\mathcal{T}_\varphi$
is finite or $\mathcal{T}_{\neg\varphi}$ is finite.

As in \cite{at} we permit extensions of $e$-minimal /
$e$-categorical families by their accumulation points and these
extensions will be also called $e$-minimal / $e$-categorical.

\medskip
{\bf Theorem 1.4} \cite{at}. {\em A family $\mathcal{T}$ is
$e$-minimal if and only if it is $e$-categorical.}

\medskip
{\bf Proposition 1.5} \cite{at}. {\em Any $E$-closed family
$\mathcal{T}$ with finite $e$-${\rm Sp}(\mathcal{T})>0$ is
represented as a disjoint union of $e$-categorical families
$\mathcal{T}_1,\ldots,\mathcal{T}_n$.}

\medskip
Proof. Let $e$-${\rm Sp}(\mathcal{T})=n$ and $T_1,\ldots,T_n$ be
accumulation points for $\mathcal{T}$ witnessing that equality.
Now we consider pairwise inconsistent formulas $\varphi_i\in T_i$
separating $T_i$ from $T_j$, $j\ne i$, i.e., with
$\neg\varphi_i\in T_j$. By Proposition 1.1 each family
$\mathcal{T}_i=\mathcal{T}_{\varphi_i}$ is infinite, with unique
accumulation point $T_i$, and thus $\mathcal{T}_i$ is
$e$-categorical. Besides, the families $\mathcal{T}_i$ are
disjoint by the choice of $\varphi_i$, and
$\mathcal{T}'=\mathcal{T}\setminus\left(\bigcup\limits_{i=1}^n\mathcal{T}_i\right)$
does not have accumulation points. Therefore
$\mathcal{T}'\cup\mathcal{T}_1$ is $e$-categorical, too. Thus,
$\mathcal{T}'\cup\mathcal{T}_1,\mathcal{T}_2,\ldots,\mathcal{T}_n$
is the required partition of $\mathcal{T}$ on $e$-categorical
families.~$\Box$

\medskip
{\bf Theorem 1.6} \cite{at}. {\em A family $\mathcal{T}$ of
theories contains an approximating subfamily if and only if
$\mathcal{T}$ is infinite.}

\medskip
Proof. Since any approximating family is infinite then, having an
approximating subfamily, $\mathcal{T}$ is infinite.

Conversely, let $\mathcal{T}$ be infinite. Firstly, we assume that
the language $\Sigma=\Sigma(\mathcal{T})$ of $\mathcal{T}$ is at
most countable. We enumerate all $\Sigma$-sentences: $\varphi_n$,
$n\in\omega$, and construct an accumulation point for
$\mathcal{T}$ by induction. Since $\mathcal{T}_{\varphi_0}$ or
$\mathcal{T}_{\neg\varphi_0}$ is infinite  we can choose
$\psi_0=\varphi^{\delta}_0$ with infinite
$\mathcal{T}_{\varphi^\delta_0}$, $\delta\in\{0,1\}$. If $\psi_n$
is already defined, with infinite $\mathcal{T}_{\psi_n}$, then we
choose $\psi_{n+1}=\psi_n\wedge\varphi^\delta_{n+1}$, with
$\delta\in\{0,1\}$, such that $\mathcal{T}_{\psi_{n+1}}$ is
infinite. Finally, the set $\{\psi_n\mid n\in\omega\}$ forces a
complete theory $T$ being an accumulation point both for
$\mathcal{T}$ and for each $\mathcal{T}_{\psi_n}$. Thus,
$\mathcal{T}\setminus\{T\}$ is a required approximating family.

If $\Sigma$ is uncountable we find an accumulation point $T_0$ for
infinite $\mathcal{T}\upharpoonright\Sigma_0$, where $\Sigma_0$ is
a countable sublanguage of $\Sigma$. Now we extend $T_0$ till a
complete $\Sigma$-theory $T$ adding $\Sigma$-sentences $\chi$ such
that $\mathcal{T}_\chi$ are infinite. Again
$\mathcal{T}\setminus\{T\}$ is a required approximating
family.~$\Box$

\section{Ranks and $e$-spectra}

Starting with $e$-categorical, i.e., $e$-minimal families of
theories we define the {\em rank} ${\rm RS}(\cdot)$ for the
families of theories, similar to Morley rank \cite{Morley}, and a
hierarchy with respect to these ranks in the following way.

For the empty family $\mathcal{T}$ we put the rank ${\rm
RS}(\mathcal{T})=-1$, for finite nonempty families $\mathcal{T}$
we put ${\rm RS}(\mathcal{T})=0$, and for infinite families
$\mathcal{T}$~--- ${\rm RS}(\mathcal{T})\geq 1$.

For a family $\mathcal{T}$ and an ordinal $\alpha=\beta+1$ we put
${\rm RS}(\mathcal{T})\geq\alpha$ if there are pairwise
inconsistent $\Sigma(\mathcal{T})$-sentences $\varphi_n$,
$n\in\omega$, such that ${\rm
RS}(\mathcal{T}_{\varphi_n})\geq\beta$, $n\in\omega$.

If $\alpha$ is a limit ordinal then ${\rm
RS}(\mathcal{T})\geq\alpha$ if ${\rm RS}(\mathcal{T})\geq\beta$
for any $\beta<\alpha$.

We set ${\rm RS}(\mathcal{T})=\alpha$ if ${\rm
RS}(\mathcal{T})\geq\alpha$ and ${\rm
RS}(\mathcal{T})\not\geq\alpha+1$.

If ${\rm RS}(\mathcal{T})\geq\alpha$ for any $\alpha$, we put
${\rm RS}(\mathcal{T})=\infty$.

A family $\mathcal{T}$ is called {\em $e$-totally transcendental},
or {\em totally transcendental}, if ${\rm RS}(\mathcal{T})$ is an
ordinal.

\medskip
Clearly, there are many totally transcendental families. At the
same time, the following example shows that there are families
which are not totally transcendental.

\medskip
{\bf Example 2.1.} Let $\mathcal{T}$ be a family of all theories,
with infinite models, in the language $\Sigma=\{Q_n\mid
n\in\omega\}$ of unary predicates such that any $Q_n$ is either
empty or complete, each $T\in\mathcal{T}$ has infinitely and
co-infinitely empty predicates, and each infinite a and
co-infinite $\Sigma_0\subset\Sigma$ has a theory $T\in\mathcal{T}$
such that $Q_n=\emptyset$ for $T$ if and only if $Q_n\in\Sigma_0$.

Since each $\Sigma$-sentence $\varphi$ is reduced to a description
of finitely many $Q_n$ that some of them are (non)empty, we always
can divide $\mathcal{T}_\varphi$ into infinitely many disjoint
parts with respect to some formulas. It implies that ${\rm
RS}(\mathcal{T})>\alpha$ for any ordinal $\alpha$, i.e.,
$\mathcal{T}$ is not totally transcendental.~$\Box$

\medskip
By the definition, since there are ${\rm
max}\{|\Sigma(\mathcal{T})|,\omega\}$
$\Sigma(\mathcal{T})$-sentences, so if ${\rm
RS}(\mathcal{T})<\infty$ then $|{\rm RS}(\mathcal{T})|\leq{\rm
max}\{|\Sigma(\mathcal{T})|,\omega\}$.

In particular, the following proposition holds.

\medskip
{\bf Proposition 2.2.} {\em If $|\Sigma(\mathcal{T})|\leq\omega$
then either $|{\rm RS}(\mathcal{T})|\leq\omega$ or $\mathcal{T}$
is not $e$-totally transcendental.}

\medskip
If $\mathcal{T}$ is totally transcendental, with ${\rm
RS}(\mathcal{T})=\alpha\geq 0$, we define the {\em degree} ${\rm
ds}(\mathcal{T})$ of $\mathcal{T}$ as the maximal number of
pairwise inconsistent sentences $\varphi_i$ such that ${\rm
RS}(\mathcal{T}_{\varphi_i})=\alpha$.

\medskip
Clearly, if ${\rm RS}(\mathcal{T})=\alpha$ then ${\rm
ds}(\mathcal{T})\in\omega\setminus\{0\}$.

Notice also that the rank ${\rm RS}(\cdot)$ is monotone both with
respect to extensions of $\mathcal{T}$ and expansions of theories
in $\mathcal{T}$: if $\mathcal{T}_1\subseteq\mathcal{T}_2$ or
$\mathcal{T}_2$ is obtained from $\mathcal{T}_1$ by expansions of
theories in $\mathcal{T}_1$ then ${\rm RS}(\mathcal{T}_1)\leq{\rm
RS}(\mathcal{T}_2)$. Besides, if ${\rm RS}(\mathcal{T}_1)$ is an
ordinal and ${\rm RS}(\mathcal{T}_1)={\rm RS}(\mathcal{T}_2)$ then
${\rm ds}(\mathcal{T}_1)\leq{\rm ds}(\mathcal{T}_2)$.

The following proposition is obvious.

\medskip
{\bf Proposition 2.3.} {\em A family $\mathcal{T}$ is $e$-minimal
if and only if ${\rm RS}(\mathcal{T})=1$ and ${\rm
ds}(\mathcal{T})=1$.}

\medskip
Thus, we have an additional, with respect to Theorem 1.4,
characterization of $e$-categoricity in terms of ranks.

\medskip
{\bf Remark 2.4.} Clearly, if ${\rm RS}(\mathcal{T})>0$ is an
ordinal then $\mathcal{T}$ can be expanded obtaining a family
$\mathcal{T}'$ such that ${\rm RS}(\mathcal{T}')>{\rm
RS}(\mathcal{T})$. Indeed, each $e$-minimal subfamily
$\mathcal{T}_\varphi$ of $\mathcal{T}$ can be divided into
countably many infinite parts just introducing countably many new
predicate such that these predicates are either empty or complete
and for any partition of $\mathcal{T}$ into countably many
infinite parts $\mathcal{T}_i$ each part can be labelled by a
sentence that some new predicate in nonempty. This procedure
increase finite rank ${\rm RS}(\mathcal{T})$ till ${\rm
RS}(\mathcal{T})+1$. If ${\rm RS}(\mathcal{T})$ is infinite, we
increase this rank either continuing to divide $e$-minimal
$\mathcal{T}_\varphi$ and obtaining ${\rm RS}(\mathcal{T}')$, or
using similar expansions by new empty and complete predicates
preserving $e$-minimality but increasing possibilities of other
steps including limit ones and obtaining an ordinal ${\rm
RS}(\mathcal{T}')>{\rm RS}(\mathcal{T})$.

\medskip
{\bf Proposition 2.5.} {\em For any infinite family $\mathcal{T}$,
$e$-${\rm Sp}(\mathcal{T})$ is finite if and only if ${\rm
RS}(\mathcal{T})=1$. If ${\rm RS}(\mathcal{T})=1$ then $e$-${\rm
Sp}(\mathcal{T})={\rm ds}(\mathcal{T})$.}

\medskip
Proof. If $e$-${\rm Sp}(\mathcal{T})$ is finite then ${\rm
RS}(\mathcal{T})=1$ following the proof of Proposition 1.5.
Conversely, if ${\rm RS}(\mathcal{T})=1$ then $\mathcal{T}$ is
divided onto ${\rm ds}(\mathcal{T})$ disjoint $e$-minimal
subfamilies $\mathcal{T}_\varphi$, each of which has a proper
accumulation point, and all accumulation points for $\mathcal{T}$
are exhausted by these ones. Thus, $e$-${\rm Sp}(\mathcal{T})={\rm
ds}(\mathcal{T})$.~$\Box$

\medskip
{\bf Proposition 2.6.} {\em If ${\rm RS}(\mathcal{T})\geq 2$ then
$e$-${\rm Sp}(\mathcal{T})\geq\omega$.}

\medskip
Proof. Since each infinite neighbourhood $\mathcal{T}_\varphi$ has
an accumulation point $T$ containing $\varphi$ and by ${\rm
RS}(\mathcal{T})\geq 2$ there are infinitely many disjoint
infinite neighbourhoods $\mathcal{T}_\varphi$ we have infinitely
many accumulation points each of which can be counted for the
value $e$-${\rm Sp}(\mathcal{T})$. Thus, $e$-${\rm
Sp}(\mathcal{T})\geq\omega$.~$\Box$

\medskip
The following example of a family $\mathcal{T}$ with ${\rm
RS}(\mathcal{T})=2$ illustrates an existence of an accumulation
point $T$ such that $T\notin{\rm Cl}_E(\mathcal{T}_{\varphi_n})$,
$n\in\omega$, where the families $\mathcal{T}_{\varphi_n}$ divide
$\mathcal{T}$ disjointly on $e$-minimal parts.

\medskip
{\bf Example 2.7.} Let $\mathcal{T}$ consists of theories $T^k_n$
with infinite models, in the language $\Sigma=\{Q^k\mid k\in\omega
\}\cup\{Q^k_n\mid k,n\in\omega\}$ of unary predicates each of
which is either empty or complete, where $Q^m$ and $Q^m_s$ are
empty for $T^k_n$, if $m\ne k$, $Q^k$ is complete for $T^k_n$, and
$Q^k_m$ is complete for $T^k_n$ if and only if $m\geq n$. The
neighbourhoods $\mathcal{T}_{\varphi_k}$, where $\varphi_k$
witnesses that $Q^k$ is complete, divide $\mathcal{T}$ onto
countably many parts, each of which is $e$-minimal. At the same
time $\mathcal{T}$ has an additional accumulation point $T$, whose
all unary predicates are empty.~$\Box$

\medskip
More generally, if there are infinitely many infinite (in
particular, $e$-minimal) families $\mathcal{T}_{\varphi_k}$ with
pairwise inconsistent sentences $\varphi_k$, $k\in\omega$, then
$\mathcal{T}$ has an accumulation point $T$ containing
$\{\neg\varphi_k\mid k\in\omega\}$. Indeed, each family
$\mathcal{T}_{\psi_n}$, where
$\psi_n=\neg\varphi_0\wedge\ldots\wedge\neg\varphi_n$, is
infinite. Therefore we can construct $T$ repeating the arguments
for Theorem 1.6.

Thus, we have the following:

\medskip
{\bf Proposition 2.8.} {\em If there are infinite families
$\mathcal{T}_{\varphi_k}$ with pairwise inconsistent sentences
$\varphi_k$, $k\in\omega$, {\rm (}witnessing ${\rm
RS}(\mathcal{T})\geq 2${\rm )} then there is an accumulation point
$T$ for $\mathcal{T}$ which is not an accumulation point for any
$\mathcal{T}_{\varphi_k}$, $k\in\omega$.}

\medskip
The following modification of Example 2.7 shows that, having ${\rm
RS}(\mathcal{T})=2$, the number of accumulation points in
$\bigcap\limits_{k\in\omega}\mathcal{T}_{\neg\varphi_k}$ can vary
from $1$ to $\omega$.

\medskip
{\bf Example 2.9.} Obtaining $n\in\omega$ additional accumulation
points it suffices take the family $\mathcal{T}$ in Example 2.7
and to mark exactly one theory in each $\mathcal{T}_{\varphi_k}$
by some new complete predicate $R_i$ such that new accumulation
point has exactly one complete predicate $R_i$. Clearly, we can
mark $n$ disjoint sequences of theories producing $n$ new
accumulation points. And it is possible to continue this process
obtaining a family $\mathcal{T}'$ with $\omega$ accumulation
points. This process preserves $e$-minimality for
$\mathcal{T}'_{\varphi_k}$ and gives the values ${\rm
RS}(\mathcal{T}')=2$ and ${\rm ds}(\mathcal{T}')=1$.~$\Box$

\medskip
It is easy to see that Example 2.9 can be naturally modified for
an arbitrarily large language, by additional complete and empty
predicates $R_i$ such that exactly one $R_i$ is complete for a
chosen theory, producing a family $\mathcal{T}$ with ${\rm
RS}(\mathcal{T})=2$, ${\rm ds}(\mathcal{T})=1$ and $e$-${\rm
Sp}(\mathcal{T})$ equals a chosen cardinality $\lambda>\omega$.

\medskip
{\bf Theorem 2.10.} {\em For any family $\mathcal{T}$, ${\rm
RS}(\mathcal{T})={\rm RS}({\rm Cl}_E(\mathcal{T}))$, and if
$\mathcal{T}$ is nonempty and $e$-totally transcendental then
${\rm ds}(\mathcal{T})={\rm ds}({\rm Cl}_E(\mathcal{T}))$.}

\medskip
Proof. At first we argue to show that ${\rm RS}(\mathcal{T})={\rm
RS}({\rm Cl}_E(\mathcal{T}))$. Since ${\rm
RS}(\mathcal{T}_1)\leq{\rm RS}(\mathcal{T}_2)$ for
$\mathcal{T}_1\subseteq\mathcal{T}_2$, and
$\mathcal{T}\subseteq{\rm Cl}_E(\mathcal{T})$, we have ${\rm
RS}(\mathcal{T})\leq{\rm RS}({\rm Cl}_E(\mathcal{T}))$. Now we
will prove the inequality \begin{equation}\label{seq1}{\rm
RS}({\rm Cl}_E(\mathcal{T}))\leq{\rm
RS}(\mathcal{T})\end{equation} by induction. If $\mathcal{T}$ is
finite then ${\rm Cl}_E(\mathcal{T})=\mathcal{T}$ and the
inequality (\ref{seq1}) is obvious. If ${\rm RS}(\mathcal{T})=1$
then by Propositions 2.3 and 2.5, $\mathcal{T}$ is a finite (with
${\rm ds}(\mathcal{T})$ parts) disjoint union of $e$-minimal,
i.e., $e$-categorical families $\mathcal{T}_\varphi$ such that
${\rm Cl}_E(\mathcal{T})=\bigcup\limits_\varphi{\rm
Cl}_E(\mathcal{T}_\varphi)$ and $|{\rm
Cl}_E(\mathcal{T})\setminus\mathcal{T}|\leq{\rm
ds}(\mathcal{T})<\omega$. Then ${\rm Cl}_E(\mathcal{T})$ is a
finite disjoint union of ${\rm ds}(\mathcal{T})$ $e$-minimal
families ${\rm Cl}_E(\mathcal{T}_\varphi)$ producing the
inequality (\ref{seq1}), with ${\rm ds}(\mathcal{T})={\rm ds}({\rm
Cl}_E(\mathcal{T}))$.

If ${\rm RS}({\rm Cl}_E(\mathcal{T}))\geq\alpha$ for a limit
ordinal $\alpha$ then ${\rm RS}(\mathcal{T})\geq\alpha$ by
induction. So it suffices to observe ${\rm
RS}(\mathcal{T})\geq\alpha+1$ if ${\rm RS}({\rm
Cl}_E(\mathcal{T}))\geq\alpha+1$. But if the latter inequality is
witnessed by some sentences $\varphi_n$, $n\in\omega$, with ${\rm
RS}({\rm Cl}_E(\mathcal{T}_{\varphi_n}))\geq\alpha$ then by
induction ${\rm RS}(\mathcal{T}_{\varphi_n})\geq\alpha$, with
${\rm Cl}_E(\mathcal{T}_{\varphi_n})={\rm
Cl}_E(\mathcal{T})_{\varphi_n}$. Therefore, ${\rm
RS}(\mathcal{T})\geq\alpha+1$ witnessed by the same sentences
$\varphi_n$.

Thus, ${\rm RS}(\mathcal{T})={\rm RS}({\rm Cl}_E(\mathcal{T}))$.

The condition ${\rm ds}(\mathcal{T})=k\Leftrightarrow{\rm ds}({\rm
Cl}_E(\mathcal{T}))=k$ follows again by the equality ${\rm
Cl}_E(\mathcal{T}_{\varphi_n})={\rm
Cl}_E(\mathcal{T})_{\varphi_n}$, where the $E$-closures of
disjoint neighbourhoods $\mathcal{T}_{\varphi_n}$, with ${\rm
RS}(\mathcal{T}_{\varphi_n})={\rm RS}(\mathcal{T})$, ${\rm
ds}(\mathcal{T}_{\varphi_n})=1$, exhaust ${\rm
Cl}_E(\mathcal{T})$.~$\Box$

\medskip
Notice that Example 2.7 can be naturally generalized in a
countable language of unary predicates producing a family
$\mathcal{T}$ with given countable ordinal $\alpha={\rm
RS}(\mathcal{T})$ and given positive natural number $n={\rm
ds}(\mathcal{T})$. Thus, the hierarchy of families $\mathcal{T}$,
in countable languages, with respect to pairs pairs
$(\alpha,n)=({\rm RS}(\mathcal{T}),{\rm ds}(\mathcal{T}))$ can be
realized.

If the language $\Sigma$ is uncountable we can continue the
process increasing ${\rm RS}(\mathcal{T})$ to uncountable ordinals
with an upper bound $|\Sigma|$, since this bound equals the
cardinality of the set of all $\Sigma$-sentences $\varphi$,
defining $\mathcal{T}_\varphi$.

Therefore the following proposition holds.

\medskip
{\bf Proposition 2.11.} {\em For any ordinal $\alpha$ and a
natural number $n\in\omega\setminus\{0\}$ there is a family $T$
such that $({\rm RS}(\mathcal{T}),{\rm
ds}(\mathcal{T}))=(\alpha,n)$.}

\medskip
Having a hierarchy with $({\rm RS}(\mathcal{T}),{\rm
ds}(\mathcal{T}))=(\alpha,n)$ and Proposition 2.5 for $({\rm
RS}(\mathcal{T}),{\rm ds}(\mathcal{T}))=(1,n)$, it is natural to
characterize these values $(\alpha,n)$ for $\alpha\geq 2$.

\medskip
{\bf Definition.} A family $\mathcal{T}$, with infinitely many
accumulation points, is called {\em $a$-minimal} if for any
sentence $\varphi\in\Sigma(T)$, $\mathcal{T}_\varphi$ or
$\mathcal{T}_{\neg\varphi}$ has finitely many accumulation points.

\medskip
The following theorem gives a characterization, in terms of
$a$-minimality, for ${\rm RS}(\mathcal{T})=2$. Notice that by
Theorem 2.10 it does not matter $\mathcal{T}$ is $E$-closed or
not.

\medskip
{\bf Theorem 2.12.} {\em For any family $\mathcal{T}$, ${\rm
RS}(\mathcal{T})=2$, with ${\rm ds}(\mathcal{T})=n$, if and only
if $\mathcal{T}$ is represented as a disjoint union of subfamilies
$\mathcal{T}_{\varphi_1},\ldots,\mathcal{T}_{\varphi_n}$, for some
pairwise inconsistent sentences $\varphi_1,\ldots,\varphi_n$, such
that each $\mathcal{T}_{\varphi_i}$ is $a$-minimal.}

\medskip
Proof. Let ${\rm RS}(\mathcal{T})=2$ and ${\rm
ds}(\mathcal{T})=n$. By the definition $\mathcal{T}$ is
represented as a disjoint union of subfamilies
$\mathcal{T}_{\varphi_1},\ldots,\mathcal{T}_{\varphi_n}$, for some
sentences $\varphi_1,\ldots,\varphi_n$, such that each
$\mathcal{T}_{\varphi_i}$ satisfies ${\rm
RS}(\mathcal{T}_{\varphi_i})=2$ and ${\rm
ds}(\mathcal{T}_{\varphi_i})=1$. So it suffices to show that,
assuming ${\rm ds}(\mathcal{T})=1$, ${\rm RS}(\mathcal{T})=2$ if
and only if $\mathcal{T}$ is $a$-minimal.

Let ${\rm RS}(\mathcal{T})=2$ and ${\rm ds}(\mathcal{T})=1$.
Therefore have infinitely many accumulation points belonging to
the $E$-closures of $e$-minimal subfamilies $\mathcal{T}_\varphi$.
If $\mathcal{T}$ is not $a$-minimal then for some sentence
$\psi\in\Sigma(T)$, $\mathcal{T}_\psi$ and
$\mathcal{T}_{\neg\psi}$ have infinitely many accumulation points.
By Proposition 2.5, ${\rm RS}\left(\mathcal{T}_\varphi\right)=2$
and ${\rm RS}\left(\mathcal{T}_{\neg\varphi}\right)=2$
contradicting ${\rm ds}(\mathcal{T})=1$.

Now let $\mathcal{T}$ be $a$-minimal. Having infinitely many
accumulation points for $\mathcal{T}$ it is easy to construct
step-by-step infinitely many disjoint infinite subfamilies
$\mathcal{T}_{\psi_i}$, $i\in\omega$, with pairwise inconsistent
sentences $\psi_i$, witnessing ${\rm RS}(\mathcal{T})\geq 2$.
Moreover, since $\mathcal{T}$ is $a$-minimal it is possible to
choose $\psi_i$ such that each $\mathcal{T}_{\psi_i}$ has unique
accumulation point, i.e., by Theorem 1.4 and Proposition 2.3, it
is $e$-minimal with ${\rm RS}\left(\mathcal{T}_{\psi_i}\right)=1$
and ${\rm ds}\left(\mathcal{T}_{\psi_i}\right)=1$. And each
possibility to divide $\mathcal{T}$ by sentences witnessing ${\rm
RS}(\mathcal{T})\geq 2$ is reduced to the case above. It means
that ${\rm RS}(\mathcal{T})=2$. Since, by $a$-minimality,
$\mathcal{T}$ can not be divided, by a sentence $\chi$, to
subfamilies $\mathcal{T}_\chi$ and $\mathcal{T}_{\neg\chi}$ with
infinitely many accumulation points, ${\rm
ds}(\mathcal{T})=1$.~$\Box$

\medskip
Below we generalize the notions of $e$-minimality and
$a$-minimality for arbitrary nonempty $e$-totally transcendental
families $\mathcal{T}$ allowing to characterize step-by-step
families of ranks $\alpha$ starting with $\alpha\in\{0,1,2\}$.

\medskip
{\bf Definition.} Let $\alpha$ be an ordinal. A family
$\mathcal{T}$ of rank $\alpha$ is called {\em $\alpha$-minimal} if
for any sentence $\varphi\in\Sigma(T)$, ${\rm
RS}(\mathcal{T}_\varphi)<\alpha$ or ${\rm
RS}(\mathcal{T}_{\neg\varphi})<\alpha$.

\medskip
By the definition and in view of Proposition 2.3 and Theorem 2.12
we have:

\medskip
{\bf Proposition 2.13.} {\em $(1)$ A family $\mathcal{T}$ is
$0$-minimal if and only if $\mathcal{T}$ is a singleton.

$(2)$ A family $\mathcal{T}$ is $1$-minimal if and only if
$\mathcal{T}$ is $e$-minimal.

$(3)$ A family $\mathcal{T}$ is $2$-minimal if and only if
$\mathcal{T}$ is $a$-minimal.

$(4)$ For any ordinal $\alpha$ a family $\mathcal{T}$ is
$\alpha$-minimal if and only if ${\rm RS}(\mathcal{T})=\alpha$ and
${\rm ds}(\mathcal{T})=1$. }

\medskip
In view of Proposition 2.13 the following assertion obviously
generalizes Theorem 2.12.

\medskip
{\bf Proposition 2.14.} {\em For any family $\mathcal{T}$, ${\rm
RS}(\mathcal{T})=\alpha$, with ${\rm ds}(\mathcal{T})=n$, if and
only if $\mathcal{T}$ is represented as a disjoint union of
subfamilies
$\mathcal{T}_{\varphi_1},\ldots,\mathcal{T}_{\varphi_n}$, for some
pairwise inconsistent sentences $\varphi_1,\ldots,\varphi_n$, such
that each $\mathcal{T}_{\varphi_i}$ is $\alpha$-minimal.}

\section{Boolean algebras and CB-ranks}

\medskip
Similarly \cite{Morley}, for a nonempty family $\mathcal{T}$, we
denote by $\mathcal{B}(\mathcal{T})$ the Boolean algebra
consisting of all subfamilies $\mathcal{T}_{\varphi}$, where
$\varphi$ are sentences in the language $\Sigma(\mathcal{T})$.

Following \cite{Morley} we observe that $\mathcal{B}(\mathcal{T})$
is superatomic \cite{Day, BA} for every $e$-totally transcendental
$\mathcal{T}$, with well-ordered chains. And vice versa, having
superatomic $\mathcal{B}(\mathcal{T})$ we step-by-step define
ordinals ${\rm RS}(\mathcal{T}_{\varphi})$ for
$\mathcal{T}_{\varphi}$ implying that $\mathcal{T}$ is $e$-totally
transcendental. Thus, the following theorem holds.

\medskip
{\bf Theorem 3.1.} {\em A nonempty family $\mathcal{T}$ is
$e$-totally transcendental if and only if the Boolean algebra
$\mathcal{B}(\mathcal{T})$ is superatomic.}

\medskip
In particular, for an infinite family $\mathcal{T}$, the start of
the process, producing an ordinal ${\rm RS}(\mathcal{T})$, should
be bases on $e$-minimal families $\mathcal{T}_{\varphi}$, i.e., if
each infinite $\mathcal{T}_{\varphi}$ is definably divided into
two infinite parts $\mathcal{T}_{\varphi\wedge\psi}$ and
$\mathcal{T}_{\varphi\wedge\neg\psi}$, then
$\mathcal{T}_{\varphi}$, and, in particular,
$\mathcal{T}=\mathcal{T}_{\forall x(x\approx x)}$, has ${\rm
RS}(\mathcal{T}_{\varphi})=\infty$.

Thus we have the following

\medskip
{\bf Proposition 3.2.} {\em If an infinite family $\mathcal{T}$
does not have $e$-minimal subfamilies $\mathcal{T}_\varphi$ then
$\mathcal{T}$ is not $e$-totally transcendental.}

\medskip
{\bf Remark 3.3.} By the definition of the rank, for any family
$\mathcal{T}$ represented as a union
$\mathcal{T}_1\cup\mathcal{T}_2$ we have ${\rm
RS}(\mathcal{T})={\rm max}\{{\rm RS}(\mathcal{T}_1),{\rm
RS}(\mathcal{T}_2)\}$ since each step for ${\rm RS}(\mathcal{T})$
uses infinitely many theories in $\mathcal{T}_1$ or
$\mathcal{T}_2$ dividing some neighbourhoods
$(\mathcal{T}_i)_{\varphi}$ into infinitely many disjoint parts.
At the same time, ${\rm ds}(\mathcal{T})$ can vary from ${\rm
max}\{{\rm ds}(\mathcal{T}_1),{\rm ds}(\mathcal{T}_2)\}$ till
${\rm ds}(\mathcal{T}_1)+{\rm ds}(\mathcal{T}_2)$ depending on
$\mathcal{T}_1$ and $\mathcal{T}_2$.~$\Box$

\medskip
Recall the definition of the Cantor--Bendixson rank. It is defined
on the elements of a topological space $X$ by induction: ${\rm
CB}_X(p)\geq 0$ for all $p\in X$; ${\rm CB}_X(p)\geq\alpha$ if and
only if for any $\beta<\alpha$, $p$ is an accumulation point of
the points of ${\rm CB}_X$-rank at least $\beta$. ${\rm
CB}_X(p)=\alpha$ if and only if both ${\rm CB}_X(p)\geq\alpha$ and
${\rm CB}_X(p)\ngeq\alpha+1$ hold; if such an ordinal $\alpha$
does not exist then ${\rm CB}_X(p)=\infty$. Isolated points of $X$
are precisely those having rank $0$, points of rank $1$ are those
which are isolated in the subspace of all non-isolated points, and
so on. For a non-empty $C\subseteq X$  we define ${\rm
CB}_X(C)=\sup\{{\rm CB}_X(p)\mid p\in C\}$; in this way ${\rm
CB}_X(X)$ is defined and ${\rm CB}_X(\{p\})={\rm CB}_X(p)$ holds.
If $X$ is compact and $C$ is closed in $X$ then the sup is
achieved: ${\rm CB}_X(C)$ is the maximum value of ${\rm CB}_X(p)$
for $p\in C$; there are  finitely many points of maximum rank in
$C$ and the number of such points is the \emph{${\rm
CB}_X$-degree} of $C$, denoted by $n_X(C)$.

If $X$ is countable and compact then ${\rm CB}_X(X)$ is a
countable ordinal and every  closed subset has ordinal-valued rank
and finite ${\rm CB}_X$-degree $n_X(X)\in\omega\setminus\{0\}$.

For any ordinal $\alpha$ the set $\{p\in X\mid{\rm
CB}_X(p)\geq\alpha\}$ is called the $\alpha$-th {\em ${\rm
CB}$-derivative} $X_\alpha$ of $X$.

Elements $p\in X$ with ${\rm CB}_X(p)=\infty$ form the {\em
perfect kernel} $X_\infty$ of $X$.

Clearly, $X_\alpha\supseteq X_{\alpha+1}$, $\alpha\in{\rm Ord}$,
and $X_\infty=\bigcap\limits_{\alpha\in{\rm Ord}}X_\alpha$.

Similarly, for a nontrivial superatomic Boolean algebra
$\mathcal{A}$ the characteristics ${\rm CB}_\mathcal{A}(A)$,
$n_\mathcal{A}(A)$, and ${\rm CB}_\mathcal{A}(p)$, for $p\in A$,
are defined \cite{BA} starting with atomic elements being isolated
points. Following \cite{BA}, ${\rm CB}_\mathcal{A}(A)$ and
$n_\mathcal{A}(A)$ are called the {\em Cantor--Bendixson
invariants}, or {\em ${\rm CB}$-invariants} of $\mathcal{A}$.

Recall that by \cite[Lemma 17.9]{BA}, ${\rm
CB}_\mathcal{A}(A)<|A|^+$ for any infinite $\mathcal{A}$, and the
following theorem holds.

\medskip
{\bf Theorem 3.4} \cite[Theorem 17.11]{BA}. {\em Countable
superatomic Boolean algebras are isomorphic if and only if they
have the same ${\rm CB}$-invariants.}

\medskip
In view of Theorem 3.1 any $e$-totally transcendental family
$\mathcal{T}$ defines a superatomic Boolean algebra
$\mathcal{B}(\mathcal{T})$, and it is easy to observe step-by-step
that ${\rm RS}(\mathcal{T})={\rm
CB}_{\mathcal{B}(\mathcal{T})}(B(\mathcal{T}))$, ${\rm
ds}(\mathcal{T})=n_{\mathcal{B}(\mathcal{T})}(B(\mathcal{T}))$,
i.e., the pair $({\rm RS}(\mathcal{T}),{\rm ds}(\mathcal{T}))$
consists of ${\rm CB}$-invariants for $\mathcal{B}(\mathcal{T})$.

In particular, by Theorem 3.4, for any countable $e$-totally
transcendental family $\mathcal{T}$, $\mathcal{B}(\mathcal{T})$ is
uniquely defined, up to isomorphism, by the pair $({\rm
RS}(\mathcal{T}),{\rm ds}(\mathcal{T}))$ of ${\rm CB}$-invariants.

By the definition for any $e$-totally transcendental family
$\mathcal{T}$ each theory $T\in\mathcal{T}$ obtains the ${\rm
CB}$-rank ${\rm CB}_\mathcal{T}(T)$ starting with
$\mathcal{T}$-isolated points $T_0$, of ${\rm
CB}_\mathcal{T}(T_0)=0$. We will denote the values ${\rm
CB}_\mathcal{T}(T)$ by ${\rm RS}_\mathcal{T}(T)$ as the rank for
the point $T$ in the topological space on $\mathcal{T}$ which is
defined with respect to $\Sigma(\mathcal{T})$-sentences.

\medskip
{\bf Remark 3.5.} By the definition we have ${\rm
RS}_\mathcal{T}(T)\geq 1$ if and only if $T$ is an accumulation
point for $\mathcal{T}$, ${\rm RS}_\mathcal{T}(T)\geq 2$ if and
only if $T$ is an accumulation point for the subfamily of ${\rm
Cl}_E(\mathcal{T})$ consisting of all its accumulation points,
etc. Additionally, by Proposition 2.14, if $\mathcal{T}$ is
$E$-closed with ${\rm RS}(\mathcal{T})=\alpha\geq 0$ then
$\mathcal{T}$ contains exactly ${\rm ds}(\mathcal{T})$ theories
$T$ such that ${\rm RS}_\mathcal{T}(T)=\alpha$. If means that
$\mathcal{T}$ is represented as a disjoint union of ${\rm
ds}(\mathcal{T})$ $\alpha$-minimal subfamilies
$\mathcal{T}_{\varphi_i}$ each of which has unique theory $T_i$
with ${\rm RS}_\mathcal{T}(T_i)=\alpha$.~$\Box$

\section{Ranks for countable languages}

Below we prove a characterization for bounds of the hierarchy of
${\rm RS}(\mathcal{T})$, for countable languages, i.e., rank
bounds for $e$-totally transcendental families.

\medskip
{\bf Proposition 4.1.} {\em If ${\rm RS}(\mathcal{T})=\infty$ then
$|{\rm Cl}_E(\mathcal{T})|\geq 2^\omega$.}

\medskip
Proof. Since ${\rm RS}(\mathcal{T})=\infty$ there is a {\em
$2$-tree} of sentences $\varphi_{\Delta}$,
$\Delta\in{\,}^{<\omega}2$, such that
$\mathcal{T}_{\varphi_{\Delta}}$ are infinite,
$\varphi_{{\Delta}\,\hat{\,}\,i}\vdash\varphi_{\Delta}$,
$i\in\{0,1\}$, and $\varphi_{{\Delta}\,\hat{\,}\,0}$,
$\varphi_{{\Delta}\,\hat{\,}\,1}$ are inconsistent. It easy to see
that for each $f\in 2^\omega$ there is an accumulation point $T_f$
for $\mathcal{T}$ containing the sentences $\varphi_{\langle
f(0),\ldots,f(n)\rangle}$, $n\in\omega$. Clearly, $T_{f_1}\ne
T_{f_2}$ for $f_1\ne f_2$. Hence, $|{\rm Cl}_E(\mathcal{T})|\geq
2^\omega$.~$\Box$

\medskip
{\bf Remark 4.2.} If the language $\Sigma(\mathcal{T})$ is at most
countable then the 2-tree of sentences $\varphi_{\Delta}$ in the
proof of Proposition 4.1 allows to form a countable subfamily
$\mathcal{T}'$ of $\mathcal{T}$ with $e$-${\rm Sp}(\mathcal{T})=
2^\omega$. For this aim it suffices to choose for $\mathcal{T}''$
some theories in $\mathcal{T}_{\varphi_{\Delta}}$ which do not
belong to some $\mathcal{T}_{\varphi_{\Delta'}}$, where $\Delta'$
is a continuation of $\Delta$. The theories $T_f$ belong to the
$E$-closure of $\mathcal{T}'$ being the union of $\mathcal{T}''$
with some at most countable subset $\mathcal{T}'''$ of
$\mathcal{T}$ such that each sentence $\varphi$ in any $T_f$ has
countable $\mathcal{T}'''_\varphi$. Thus, $e$-${\rm
Sp}(\mathcal{T})\geq 2^\omega$.

\medskip
In general case, for $|\Sigma(\mathcal{T})|\leq\omega$, both
infinite families $\mathcal{T}$ and ${\rm Cl}_E(\mathcal{T})$ are
{\em countably generated}, i.e., contain a countable
$\mathcal{T}'''$ generating both $\mathcal{T}$ and ${\rm
Cl}_E(\mathcal{T})$. Indeed, since there are countably many
$\Sigma(\mathcal{T})$-sentences $\varphi$, by Proposition 1.1 if
suffices to form $\mathcal{T}'''$ by all finite
$\mathcal{T}$-definable families $\mathcal{T}_\varphi$, and by
arbitrary countable subfamilies of $\mathcal{T}_\varphi$, if
$\mathcal{T}_\varphi$ is infinite.

Thus, for any at most countable language $\Sigma$, if $|{\rm
Cl}_E(\mathcal{T})|=2^\omega$ then $e$-${\rm
Sp}(\mathcal{T})=2^\omega$. Since $e$-${\rm
Sp}(\mathcal{T})=2^\omega$ obviously implies $|{\rm
Cl}_E(\mathcal{T})|=2^\omega$ for
$|\Sigma(\mathcal{T})|\leq\omega$, we have the following:

\medskip
{\bf Proposition 4.3.} {\em If $|\Sigma(\mathcal{T})|\leq\omega$
then $|{\rm Cl}_E(\mathcal{T})|=2^\omega$ if and only if $e$-${\rm
Sp}(\mathcal{T})=2^\omega$.}

\medskip
{\bf Proposition 4.4.} {\em If $|\Sigma(\mathcal{T})|\leq\omega$
and $|{\rm Cl}_E(\mathcal{T})|=2^\omega$ then ${\rm
RS}(\mathcal{T})=\infty$.}

\medskip
Proof. By Theorem 2.10 it suffices to assume that $\mathcal{T}$ is
$E$-closed such that $|\mathcal{T}|=2^\omega$.

At first we note that there is a sentence $\varphi$ such that
$|\mathcal{T}_\varphi|=2^\omega$ and
$|\mathcal{T}_{\neg\varphi}|=2^\omega$. Indeed, assuming that
$\varphi$ does not exist we can enumerate all
$\Sigma(\mathcal{T})$ sentences: $\varphi_n$, $n\in\omega$, and
form a sequence $\psi_n$ such that $\psi_0=\varphi_0$,
$\psi_{n+1}=\psi_n\wedge\varphi^\delta_{n+1}$, $\delta\in\{0,1\}$,
with $|\mathcal{T}_{\psi_{n+1}}|=2^\omega$. Thus,
$\left|\bigcap\limits_{n=0}^\infty\mathcal{T}_{\psi_{n}}\right|=2^\omega$
contradicting the condition that $\{\psi_n\mid n\in\omega\}$
forces a complete theory.

Repeating the arguments we construct a $2$-tree ${\rm Tr}$ of
sentences $\varphi_{\Delta}$, $\Delta\in{\,}^{<\omega}2$, as in
the proof of Proposition 4.1 such that each $\varphi_{\Delta}$
satisfies $|\mathcal{T}_{\varphi_{\Delta}}|=2^\omega$.

Now the sentences in the $2$-tree ${\rm Tr}$ witness that
$\mathcal{T}$ in not $e$-totally transcendental. Indeed,
$\mathcal{T}_{\varphi_0}$, $\mathcal{T}_{\varphi_{10}}$,
$\mathcal{T}_{\varphi_{110}},\ldots$ are disjoint families each of
which has continuum many theories. Each family
$\mathcal{T}_{\varphi_{1\ldots 10}}$ contains again infinitely
many disjoint subfamilies $\mathcal{T}_{\varphi_{1\ldots 101\ldots
10}}$ each of which has continuum many theories. Continuing the
process we observe that each $\mathcal{T}_{\varphi_{1\ldots 10}}$,
and, thus, $\mathcal{T}$ have the ranks equal to $\infty$.~$\Box$

\medskip
Collecting Propositions 4.1, 4.3, and 4.4 we obtain:

\medskip
{\bf Theorem 4.5.} {\em For any family $\mathcal{T}$ with
$|\Sigma(\mathcal{T})|\leq\omega$ the following conditions are
equivalent:

$(1)$ $|{\rm Cl}_E(\mathcal{T})|=2^\omega$;

$(2)$ $e$-${\rm Sp}(\mathcal{T})=2^\omega$;

$(3)$ ${\rm RS}(\mathcal{T})=\infty$.}

\medskip
{\bf Remark 4.6.} Having characterizations for $e$-totally
transcendental families $\mathcal{T}$ of theories by Theorem 4.5
we observe that both theories $T$ in $e$-totally transcendental
$\mathcal{T}$ can be not totally transcendental themselves,
containing, for instance, countably many independent unary
predicates, and totally transcendental theories, with either empty
or complete predicates $Q$, as in Example 2.7, can form families
$\mathcal{T}$ which are not $e$-totally transcendental, just
dividing $\mathcal{T}$ by sentences describing that the predicates
$Q$ are empty or complete. Thus, the notions of totally
transcendental theories and $e$-totally transcendental families do
not correlate in general case.

\medskip
{\bf Remark 4.7.} Examples in \cite{lut} show that families
$\mathcal{T}$ with $|\Sigma(\mathcal{T})|\leq\omega$ and $|{\rm
Cl}_E(\mathcal{T})|=2^\omega$ can (do not) have least generating
sets. Moreover, modifications of this examples can produce
families of theories with proper derivatives for arbitrary
ordinals $\alpha$. Therefore the perfect kernel for $\mathcal{T}$
can be formed on some derivative step $\alpha$. Thus, for any
ordinal $\alpha>0$ there is a family $\mathcal{T}$ such that
$\mathcal{T}_\alpha=\mathcal{T}_\infty$ whereas
$\mathcal{T}_\beta\ne\mathcal{T}_\infty$ for $\beta<\alpha$.

\medskip
{\bf Remark 4.8.} Notice that Theorem 4.5 does not hold for
$|\Sigma(\mathcal{T})|>\omega$, in general case. Indeed, language
uniform theories \cite{lut} can have both big cardinalities for
languages, big cardinalities for $\mathcal{T}$ and small
cardinalities for $e$-spectra. For instance, taking a family
$\mathcal{T}=\{T_i\mid i\in I\}$ in a language $\Sigma$ of unary
predicates $Q_i$, $i\in I$, $|I|=\lambda>\omega$, such that $T_i$
has complete predicate $Q_i$ and empty predicates $Q_j$, $j\ne i$,
we have $|{\rm Cl}_E(\mathcal{T})|=\lambda$ with ${\rm
Cl}_E(\mathcal{T})=\mathcal{T}\cup\{T_\infty\}$, where $T_\infty$
has only empty predicates, whereas $e$-${\rm Sp}(\mathcal{T})=1$,
that witnessed by $T_\infty$. Besides, $\mathcal{T}$ is
$e$-minimal, i.e., ${\rm RS}(\mathcal{T})=1$ and ${\rm
ds}(\mathcal{T})=1$. In particular, for $\lambda=2^\omega$, we
have $|{\rm Cl}_E(\mathcal{T})|=2^\omega$, $e$-${\rm
Sp}(\mathcal{T})=1$, and ${\rm RS}(\mathcal{T})=1$ refuting
Theorem 4.5 for $|\Sigma(\mathcal{T})|=2^\omega$.

Additionally, the family $\mathcal{T}$ can be expanded by unary
disjoint predicates $Q'_j$, $j\in J$, $|J|\geq 2$, such that each
$T_i$ is extended to $T_{ij}$ obtaining complete $Q'_j$ and empty
$Q'_k$ for $k\ne j$. The families $\mathcal{T}_j=\{T_{ij}\mid i\in
I\}$ stay $e$-minimal, producing unique accumulation points,
whereas we have for $\mathcal{T}'=\{T_{ij}\mid i\in I, j\in J\}$:

1) if $J$ is finite then ${\rm
Cl}_E(\mathcal{T}')=\bigcup\limits_{j\in J}{\rm
Cl}_E(\mathcal{T}_j)$ and $|{\rm Cl}_E(\mathcal{T}')|=|J|$; if $J$
is infinite then ${\rm Cl}_E(\mathcal{T}')$ consists of
$\bigcup\limits_{j\in J}{\rm Cl}_E(\mathcal{T}_j)$ and $|I|$
theories with unique nonempty $Q_i$ and all empty $Q'_j$, as well
as of unique theory $T_\infty$ with all empty predicates;
therefore $|{\rm Cl}_E(\mathcal{T}')|=|I|+|J|+1=|I|+|J|$;

2) by the previous item, $e$-${\rm Sp}(\mathcal{T})=|J|$ for
finite $J$, and $e$-${\rm Sp}(\mathcal{T})=|I|+|J|$ for infinite
$J$;

3) ${\rm RS}(\mathcal{T})=2$.

\medskip
In conclusion we formulate the following:

\medskip
{\bf Problem.} {\em Describe the rank ${\rm RS}(\cdot)$ hierarchy
for natural families of theories.}

\noindent Sobolev Institute of Mathematics, \\ 4, Acad. Koptyug
avenue, Novosibirsk, 630090, Russia; \\ Novosibirsk State
Technical
University, \\ 20, K.Marx avenue, Novosibirsk, 630073, Russia; \\
Novosibirsk State University, \\ 1, Pirogova street, Novosibirsk,
630090, Russia
\end{document}